\documentclass[12pt, leqno]{amsart}

\def\bR{{\mathbb{R}}}
\def\bR2{ {\mathbb{R}^2} }

\def\cC{{\mathcal{C}}}

\def\dRZ{$RZ$}

\def\ctC{ \tilde{\mathcal C}  }

\def\qed{\hfill $\bullet  \bullet $}

\newtheorem{theorem}            {Theorem}[section]
\newtheorem{lemma}              [theorem]{Lemma}

\newtheorem{thm}[theorem]{{\bf Theorem}}
\newtheorem{cor}[theorem]{{\bf Corollary}}

\title{Linear Matrix Inequality Representation of Sets}

\begin{document}

{\hspace{3.5in}    \today}
\bigskip
\bigskip
\bigskip
\bigskip
\bigskip
\bigskip

\maketitle

\begin{center}
\vspace{3mm}

  J. William Helton    \ \ \ \ \ \ \ \  and \ \ \ \ \ \ \      Victor Vinnikov  \ \ \ \ \ \ \

\ \ \ \ \ Mathematics Dept. UCSD \ \ \    \ \ \ \ \ \ \ \ \ \ \ \ \   Ben Gurion Univ. of the Negev

La Jolla Cal     \ \ \ \ \ \ \ \ \ \ \ \ \ \ \ \ \ \ \ \                    Beer Sheva, 84105

92093, USA           \ \ \ \ \ \ \ \        \ \ \ \ \ \ \ \  \ \ \ \ \ \ \         Israel   \ \ \ \ \ \ \ \ \ \ \ \

\ \ \ \ \ \ \ helton@math.edu     \ \ \ \ \ \          \ \ \ \ \ \ \ \  \      vinnikov@math.bgu.ac.il   \

\end{center}

\vspace{3mm}


\centerline{\bf Abstract}

This article concerns the question:
which subsets of ${\mathbb R}^m$ can be represented
with Linear Matrix Inequalities, LMIs?
This gives some perspective
on the scope and limitations of one of the most powerful
techniques commonly used in control theory.
Also before having much hope of representing engineering
problems as LMIs by automatic methods one needs
a good idea of which problems can and cannot
be represented by LMIs.
Little is currently known about such problems.
In this article we  give a necessary
condition, we call "rigid convexity", which must hold for a set ${\cC} \in {\mathbb R}^m$
in order for ${\cC}$ to have an LMI representation.
Rigid convexity is proved to be  necessary and sufficient when $m=2$.
This settles a question formally stated by Pablo Parrilo
and Berndt Sturmfels in  [PSprep].

\newpage

\section{The Problem of Representing Sets with LMI's}

Given
${\mathcal C}$ a closed   convex set in ${\mathbb R}^m$ bounded by  algebraic
hypersurfaces
$$S_{p_j}: = \{x \in {\mathbb R}^m: p_j(x) =0\},$$
with all polynomials $p_j(x) >0$ for $x$ in the interior of ${\mathcal C}$.
Which  ${\mathcal C}$ can be represented in terms
of some linear pencil $L=\{L_0, L_1, L_2, \cdots, L_m\}$ as
$${\mathcal C}^L: =
\{x = (x_1, x_2, \cdots x_m ): L_0 +  L_1  x_1 +  \cdots  +  L_m x_m
\quad \mbox{is PosSemiDef}\}?$$
We shall require by our use of the term {\bf linear}
{\bf pencil} that  $L_0, L_1, L_2, \ldots, L_m $ are symmetric real entried
matrices. Also we frequently abbreviate positive (resp. semi)
definite to PD ( resp. PSD).
A {\bf monic}  pencil is one with $L_0
= I$.
We call such a representation of a set ${\mathcal C}$ a
{\bf Linear Matrix Inequality
(LMI)  representation of ${\mathcal C}$}.

The question is formally stated by Pablo Parrilo
and Berndt Sturmfels in  [PSprep] for $m=2$,
and we resolve  the $m=2$ situation  in this paper.
At the end of the paper we speculate about generalizations and other
possible representations (Section \ref{sec:fragen}).

In this paper a  polynomial always  stands
for a polynomial with real coefficients.
Since the paper might be of interest to several audiences:
operator theory, real semi-algebraic geometry, systems engineering
and possibly partial differential equations, it is
written less tersely than is common.
Our  result here  was  announced in a survey talk 
at the conference  Mathematical Theory
of Networks and Systems  in 2002, see [Hprep].

\subsection{Motivation}

A technique which has had remarkable success in many types of
optimization problems is to convert them to
Linear Matrix Inequalities. These can be solved using semidefinite
programming algorithms
provided the number of variables $x_j$ is not extremely large. The
technique was introduced first
in the area of control in [BEFB94] and related papers and has spread
through many part of that
subject, c.f. [ESLSL00].  For a collection of applications in many
other areas see the
thesis \cite{P00} and \cite{PSprep}.

In the many applications which LMIs have found there is no systematic way
to produce LMIs
for general classes of problems.
Each area has a few special tricks which convert "lucky problems" to LMIs.
Before there is any hope of producing LMIs systematically one must have
a good idea of which types of constraint sets convert to LMIs
and which do not.
That is what this paper addresses.

After the preliminary version of this paper was circulated as a preprint,
Adrian Lewis, Pablo Parillo, and M. V. Ramana
informed us that our main result was
used in their solution,  in the positive, of a 1958 conjecture by
Peter Lax on hyperbolic polynomials arising from
the study of linear hyperbolic partial differential
equations; see \cite{LPRprep}.

\section{Solution}
\label{sec:soln}

In this section we give an algebraic statement of the solution
to the Parrilo--Sturmfels problem, see
Theorem \ref{thm:mainRZ}.
In Section \ref{sec:geom} we give an equivalent geometric statement of the solution,
see Theorem \ref{thm:Geometric}.
In Section \ref{sec:oval} we give an equivalent topological statement,
see Corollary \ref{cor:oval} and Theorem \ref{thm:oval},
which while intuitively appealing requires projective spaces
to state in full generality.
The  solution depends on results of
Victor Vinnikov \cite{V93} (see also Dubrovin \cite{Dub83}) 
and variations on it.

\subsection{A key class of polynomials}
\label{sec:RZdef}
We shall soon see that it is a special class of polynomials
we call real zero polynomials
which bound sets  amenable to LMI representation.
So now we define a {\bf real zero polynomial} ({\bf RZ} \ polynomial)
to be a polynomial in $m$ variables
satisfying for each $x\in {\mathbb R}^m$,
$$ \quad p(\mu x) = 0 \quad
\mbox{implies}\quad\mu\quad
\mbox{is real}\leqno{(RZ)}$$

A few properties are:

\noindent
{\bf RZ Properties}
\begin{enumerate}

\item{ The product $p_1 p_2 \cdots p_k$ of \dRZ  polynomials $p_1, \ldots, p_k$
is a  \dRZ \ polynomial.}
\item{ If a \dRZ \ polynomial $p$ factors as $p=p_1 p_2 \cdots p_k$, then
all factors $p_1, p_2, \ldots, p_k$
are \dRZ \ polynomials.}
%
%
\item{ The determinant
\begin{equation}
\label{eq:detL}
{\check p}(x) := \det [L_0 + L_1 x_1 + \cdots +  L_m x_m]
\end{equation}
of a pencil $L$ with a PD matrix $L_0$,  is an
  \dRZ \ polynomial.}
%
%
\item{The product of ${\check p}_j$ defined as determinants of pencils $L^j$
on ${\mathbb R}^{N_j}$, to wit
\begin{equation}
{\check p}_j(x) := \det [L_0^j + L_1^j x_1 + \cdots +  L_m^j x_m]
\end{equation}
is also the determinant of a pencil $L$, namely, of the direct sum
  $L_i: = L^1_i \oplus \cdots \oplus L^k_i$ on a single
${\mathbb R}^N$ with $N=N_1 + \dots + N_k$.
We have
\begin{eqnarray*}
\det L &=& (\det L^1) ( \det L^2) \dots (\det L^k)\\
&=& {\check p}_1 {\check p}_2 \dots {\check p}_k.
\end{eqnarray*}}
\item{
See Section \ref{sec:geom} and Section \ref{sec:oval} for a characterization
of \dRZ \ polynomials in terms of geometric and topological properties of
their zero sets.
The topological characterization and related topological properties of
\dRZ \ polynomials, see Section \ref{sec:Homprops},
are essentially based on connectivity results of Nuij \cite{Nuij68}.}
\end{enumerate}

\noindent
{\bf
Proofs:
}

1. 2. \ The proofs of the first two properties about products and
factors follows from the fact that the product $w= w_1 \cdots w_k$
of a set of numbers is zero if and only if  $w$ is zero.
%

3. \ For monic pencils
$$ {\check p}\left(\frac{x}{\lambda}\right) = \frac{1}{\lambda^N} \det [
\lambda I + L_1 x_1 + \cdots +  L_m x_m]$$
and all of the eigenvalues of the symmetric matrix
$L_1 x_1 + \cdots + L_mx_m$
are real, so we have that $\check p$ is \dRZ.
Replacing the monic condition by $L_0$ is PD is obvious.
A partial extension to the case when $L_0$ is merely PSD is provided by Lemma \ref{thm:monic}
below.

4. \ The proof is in the statement.\qed

We start the presentation of our main result with a key definition.

\subsection{Algebraic Interiors and their Degree}
A closed set ${\mathcal C}$ in ${\mathbb R}^m$  is an {\bf Algebraic Interior},
if there is a polynomial $p$ in $m$ variables
such that  ${\mathcal C}$ equals the closure of
a connected component of
$$\{x \in {\mathbb R}^m: p(x) > 0\}.$$
In other words, there is a polynomial $p$ in $m$ variables
which vanishes on the boundary of ${\mathcal C}$
and such that
$\{ x \in {\mathcal C}: p(x) > 0\}$ is connected with closure equal to ${\mathcal C}$.
(Notice that in general $p$ may vanish also at some points in the interior
of ${\mathcal C}$; for example, look at $p(x_1,x_2) = x_2^2 - x_1^2 (x_1 - 1)$.)
We denote the closure of the connected component of $x^0$ in
$\{x \in {\mathbb R}^m: p(x) > 0\}$
by ${\mathcal C}_p(x^0)$.
We often assume for normalization that $0$ is in the interior of ${\mathcal C}$,
and set ${\mathcal C}_p = {\mathcal C}_p(0)$.

\begin{lemma}
A polynomial $p$ of the lowest degree for which ${\mathcal C} = {\mathcal C}_p(x^0)$
is unique (up to a multiplication by a positive constant), and any other polynomial $q$
such that ${\mathcal C} = {\mathcal C}_q(x^0)$
is given by $q=ph$ where $h$ is an arbitrary polynomial which
is strictly positive on a dense connected subset of ${\mathcal C}$.
\end{lemma}

\noindent{\bf Proof}.
We shall be using some properties of algebraic and
semi-algebraic sets in ${\mathbb R}^m$, so many readers may
want to skip over it and go to our main results which are much
more widely understandable; our reference is
\cite{BCR98}.  We notice first that ${\mathcal C}$ is a semi-algebraic
set (since it is the closure of a connected component of a
semi-algebraic set, see \cite[Proposition 2.2.2 and Theorem
2.4.5]{BCR98}). Therefore the interior $\operatorname{int} {\mathcal
C}$ of ${\mathcal C}$ is also semi-algebraic, and so is the boundary
$\partial {\mathcal C} = {\mathcal C} \setminus \operatorname{int} {\mathcal C}$.
Notice also that ${\mathcal C}$ equals the closure of its interior.

We claim next that {\it for each $x \in \partial {\mathcal C}$, the
local dimension
$\dim \partial {\mathcal C}_x$ equals} $m-1$. On the one hand, we
have
$$
\dim \partial {\mathcal C}_x \leq \dim \partial {\mathcal C} <
\dim \operatorname{int} {\mathcal C} = m;
$$
here we have used \cite[Proposition 2.8.13 and Proposition 2.8.4]{BCR98},
and the fact that $\partial {\mathcal C} = \operatorname{clos} \operatorname{int} {\mathcal C}
\setminus \operatorname{int} {\mathcal C}$,
since ${\mathcal C}$ equals the closure of its interior.
On the other hand, let $B$ be an open ball in ${\mathbb R}^m$ around $x$;
then
$$
B \cap \partial {\mathcal C} = B \setminus \big[
\bigl(B \cap ({\mathbb R}^m \setminus {\mathcal C})\bigr) \cup
\bigl(B \cap \operatorname{int} {\mathcal C}\bigr)\big].
$$
Since ${\mathcal C}$ equals the closure of its interior,
every point of $\partial {\mathcal C}$ is an accumulation point of both
${\mathbb R}^m \setminus {\mathcal C}$ and $\operatorname{int} {\mathcal C}$;
therefore $B \cap ({\mathbb R}^m \setminus {\mathcal C})$
and $B \cap \operatorname{int} {\mathcal C}$ are disjoint
open nonempty semi-algebraic subsets of $B$. Using \cite[Lemma
4.5.2]{BCR98} we conclude that
$\dim B \cap \partial {\mathcal C} \geq m-1$,
hence $\dim \partial {\mathcal C}_x \geq m-1$.

Let now $V$ be the Zariski closure of $\partial {\mathcal C}$,
and let $V = V_1 \cup \cdots \cup V_k$ be the decomposition of $V$ into irreducible components.
We claim that {\it $\dim V_i = m-1$ for each $i$}.
Assume by contradiction that $V_1,\ldots,V_l$ have dimension $m-1$ while
$V_{l+1},\ldots,V_k$ have smaller dimension.
Then there exists $x \in \partial {\mathcal C}$
such that $x \not\in V_1,\ldots,V_l$, and consequently there exists an open ball $B$
in ${\mathbb R}^m$ around $x$ such that
$$
B \cap \partial {\mathcal C} = (B \cap \partial {\mathcal C} \cap V_{l+1}) \cup \cdots \cup
(B \cap \partial {\mathcal C} \cap V_k).
$$
By assumption each set in the union on the right hand side has dimension smaller than $m-1$,
hence it follows (by \cite[Proposition 2.8.5, I]{BCR98}) that
$\dim B \cap \partial {\mathcal C} < m-1$, a contradiction
with $\dim \partial {\mathcal C}_x = m-1$.

Suppose now that $p$ is a polynomial
of the lowest degree with ${\mathcal C} = {\mathcal C}_p(x^0)$.
Lowest degree implies that $p$ can have no
multiple irreducible factors, i.e., $p=p_1 \cdots p_s$, where $p_1$, \dots, $p_s$ are distinct
irreducible polynomials; we may assume without loss of generality that every $p_i$
is non-negative on ${\mathcal C}$. Since $p$ vanishes on $\partial {\mathcal C}$
it also vanishes on $V = V_1 \cup \cdots \cup V_k$.
We claim that for every $V_i$ there exists a $p_j$ so that $p_j$ vanishes on $V_i$:
otherwise ${\mathcal Z}(p_j) \cap V_i$ is a proper algebraic subset
of $V_i$ for every $j = 1,\ldots,s$, therefore (since $V_i$ is irreducible)
$\dim {\mathcal Z}(p_j) \cap V_i < \dim V_i$ for every $j$
and thus also $\dim {\mathcal Z}(p) \cap V_i < \dim V_i$,
a contradiction since $p$ vanishes on $V_i$.
If $p_j$ vanishes on $V_i$, it follows ( since $p_j$ is irreducible
and $\dim V_i = m-1$) that ${\mathcal Z}(p_j) = V_i$.
The fact that $p$ is a polynomial
of the lowest degree with ${\mathcal C} = {\mathcal C}_p(x^0)$
implies now (after possibly renumbering the irreducible factors of $p$)
that $p = p_1 \cdots p_k$ where ${\mathcal Z}(p_i) = V_i$ for every $i$.
Since $\dim V_i = m-1$ it follows from the real Nullstellensatz for principal ideals
\cite[Theorem 4.5.1]{BCR98} that the irreducible polynomials $p_i$ are uniquely determined
(up to a multiplication by a positive constant),
hence so is their product $p$.  This proves the uniqueness of $p$.

The rest of the lemma now follows easily. If ${\mathcal C} = {\mathcal C}_q(x^0)$,
then the polynomial $q$ vanishes on $\partial {\mathcal C}$ hence also on
$V = V_1 \cup \cdots \cup V_k$. Since $q$ vanishes on ${\mathcal Z}(p_i) = V_i$
and $\dim V_i = m-1$, the real Nullstellensatz for principal ideals
implies that $q$ is divisible by $p_i$; this holds for every $i$ hence
$q$ is divisible by $p = p_1 \cdots p_k$, i.e., $q = p h$.
It is obvious that $h$ must be strictly positive on a dense connected subset of ${\mathcal C}$.

\qed

\bigskip

We refer to $p$ as a {\bf minimal defining polynomial for ${\mathcal C}$},
and say that ${\mathcal C}$ is an
{\bf Algebraic Interior of Degree $d$},
where $d = \deg p$.

\subsection{Main Theorem on LMI Representations: Algebraic Version}

\begin{theorem}
\label{thm:mainRZ}


If a set
${\mathcal C}$ in ${\mathbb R}^m$, with $0$ in the interior of ${\mathcal C}$,
has an LMI representation, then
it is a convex algebraic interior and the
minimal defining polynomial $p$ for
${\mathcal C} $  satisfies the   \dRZ \ condition
(with $p(0)>0$).

Conversely, when $m=2$,
if p is an \dRZ \ polynomial of degree $d$ and  $p(0)>0$,
then ${\mathcal C}_p$ has a monic LMI representation
with $d \times d$ matrices.
\end{theorem}

The idea behind our organization of the proof is to
do the informative elementary arguments first, then
a geometric interpretation
in Section \ref{sec:geom}
(readers less interested in proofs might go there next),
and then
some examples, and then in
a seperate section do the converse side, which  refers heavily to
\cite{V93}.

\subsubsection{{ \it Proof that \dRZ  \ is Necessary }}
\label{sec:necessity}

We are given a pencil $L$ which represents ${\mathcal C}$, ${\mathcal C} = {\mathcal C}^L$;
since $0 \in {\mathcal C}$
we have $L_0$ is PSD.
Lemma \ref{thm:monic} below
reduces general pencils
with $L_0 $ a PSD matrix to monic pencils, so we assume that $L$ is a monic pencil.
Define a polynomial ${\check p}$ by
\begin{equation}
{\check p}(x) := \det [I + L_1 x_1 + \cdots +  L_m x_m].
\end{equation}
and note that \dRZ \ Property 3 says that ${\check p}$ satisfies \dRZ.
Now we show that
${\mathcal C}_{\check p}$ equals ${\mathcal C}^L$
(equals ${\mathcal C}$).
If $x(t)$, $t\in[0,1]$,
is any continuous path in ${\mathbb R}^m$ with $x(0)=0$ and $\check p(x(t)) > 0$
for all $t$, then $L(x(t))$ is a continuous path in real symmetric matrices with
$L(x(0))=I$ and $\det L(x(t)) > 0$ for all $t$, hence $L(x(1))$ is necessarily PD.
It follows that the connected component of $0$ in
$\{x \in {\mathbb R}^m: \check p(x) > 0\}$ is contained in ${\mathcal C}^L$
and therefore ${\mathcal C}_{\check p}$ is contained in ${\mathcal C}^L$.
Conversely, $x\in
{\mathcal C}^L$ implies $I+ \mu[L_1 x_1 + \cdots + L_m x_n]$ is PD for
$0 \le \mu <1$,
so $\check p(\mu x) > 0$ for $0 \le \mu <1$
and therefore necessarily $x\in {\mathcal C}_{\check p}$.

We conclude that ${\mathcal C}={\mathcal C}_{\check p}$ is a
convex algebraic interior. If $p$ is the minimal defining
polynomial for ${\mathcal C}$, then $p$ is a factor of $\check p$.
By \dRZ \ Property 2 the polynomial  $p$ is itself \dRZ. \qed

Notice that the same argument as in the proof shows that for a monic pencil $L$
the interior of ${\mathcal C} = {\mathcal C}^L$ is exactly
$\{x : I + x_1 L_1 + \cdots x_m L_m \mbox{ is PD}\}$,
and therefore the polynomial $\check p$ --- and thus also the minimal defining polynomial
$p$ --- is stricly positive on the interior of ${\mathcal C}$.

\vspace{3mm}

\subsubsection{ Reduction to Monic Pencils  }
We reduce the pencil representation problem to one
involving monic pencils as promised.

\begin{lemma}\label{thm:monic}
If ${\mathcal C}$ contains $0$ in its interior and if there is a
pencil $L$ such that $${\mathcal C} = {\mathcal C}^L,$$
then there is a monic pencil $\widehat L$ such that
$${\mathcal C}={\mathcal C}^{\widehat L}$$
\end{lemma}

\noindent{\bf Proof}.
Given ${\mathcal C}$ and $L$, since $0 \in {\mathcal C}$
we have $0 \in {\mathcal C}^L$, or equivalently, $L_0$
is PSD.  We need, however, to take $L_0$ invertible.  Since $0$
is in the interior of
${\mathcal C} \subset {\mathbb R}^m$,
it is also in the interior of $C^L$, so we have
$$
L_0 \ge \pm\, \varepsilon_1 L_1, \ \ldots,\
L_0 \ge \pm\, \varepsilon_m L_m
$$
for small enough $\varepsilon_1$  through $\varepsilon_m$.
Set
$Y=$ Range $L_0 \subset {\mathbb R}^N$, and set $$\widetilde L_j: =
L_j|_Y \quad j=0, 1, \ldots, m.$$

Clearly
$\widetilde L_0: Y \to Y$ is invertible and indeed PD.
We next show that Range $L_0$ contains Range $L_1$.
If $x \perp$ Range $L_0$, i.e., $L_0 x = 0$, then
$0 = x^T L_0 x \geq \pm\, \varepsilon_1 x^T L_1 x$ and hence
$x^T L_1 x = 0$. Since $L_0 + \varepsilon_1 L_1 \geq 0$ and
$x^T (L_0 + \varepsilon_1 L_1) x = 0$ it follows that
$(L_0 + \varepsilon_1 L_1) x = 0$, and since $L_0 x = 0$, we finally
conclude that $L_1 x = 0$, i.e., $x \perp$ Range $L_1$.

Likewise, Range $L_0$ contains the ranges of $L_2, \cdots, L_m$.
Consequently, $\widetilde L_j: Y \to Y$ are all symmetric and
${\mathcal C}^L = {\mathcal C}^{\widetilde L}$.

To build $\widehat L$, factor $\widetilde L_0 = B^TB$ with
$B$ invertible and set
$$\widehat L_j : = B^{-T} \widetilde L_j B^{-1} \quad j = 0, \ldots, m.$$
The resulting pencil $\widehat L$ is monic and ${\mathcal
C}^{\widehat L} = {\mathcal C}^{\widetilde L} = {\mathcal
C}^L$.\hfill \qed

\vspace{3mm}
\subsection{Shifted Real Zero Condition}

Assuming that $0 \in {\mathcal C}$ is just a
normalization. We can shift ${\mathbb R}^m$ by $x^0$ and obtain a shifted
polynomial
$$\tilde p(x): = p(x+x^0)$$
with ${\mathcal C}_p(x^0) - x^0 = {\mathcal C}_{\tilde p}(0)$.
Then the $RZ$ condition on $\tilde p$ shifts back to the
$RZ_{x^0}$ condition on $p$ at $x^0$.
Namely, for each $x \in {\mathbb R}^m$
$$  p(\mu x +x^0) = 0 \ \mbox{implies } \ \mu\
\mbox{is real.}\leqno{(RZ_{x_0})}
$$
An obvious consequence of our main theorem is:

\begin{cor}
\label{thm:shifted}
If a set
${\mathcal C}$ in ${\mathbb R}^m$, with $x^0$ in the interior of ${\mathcal C}$,
has an LMI representation, then
it is a convex algebraic interior and the
minimal defining polynomial $p$ for
${\mathcal C}$  satisfies the $RZ_{x^0}$ condition.

Conversely, when $m=2$,
if p is polynomial of degree $d$ satisfying $RZ_{x^0}$ and
$p(x^0)>0$, then ${\mathcal C}_p(x^0)$ has an LMI representation
with $d \times d$ matrices.
\end{cor}

\noindent
{\bf Remark}.
For $m=2$, ``going through'' LMI representations, i.e., invoking the second assertion
of Corollary \ref{thm:shifted} followed by the first assertion, we conclude that
\begin{quote}
{\it
If the polynomial $p$ satisfies $RZ_{\tilde x^0}$ for some point $\tilde
x^0$ and
$p(\tilde x^0) > 0$, then $p$ satisfies $RZ_{x^0}$ for every interior point $x^0$
of ${\mathcal C}_p(\tilde x^0)$.
}
\end{quote}
This statement holds in any dimension $m$ and can be proved, without any reference to
LMI representations, using the alternative topological characterization of $RZ$ polynomials.
See Section \ref{sec:Homprops}.

\noindent{\bf Proof of Corollary \ref{thm:shifted}}.
Assume that a set ${\mathcal C}$ in ${\mathbb R}^m$,
with $x^0$ in the interior of ${\mathcal C}$,
has a pencil representation $L$, ${\mathcal C}={\mathcal C}^L$.
Shift $L$ by $x^0$ to obtain a pencil
$\widetilde L$, that is,
$$\widetilde L_0 + \widetilde L_1 x_1
+ \dots + \widetilde L_m x_m = L_0 + L_1(x_1+x^0_1) + \dots
+ L_m (x_m+x^0_m).$$
Then ${\mathcal C}^L - x^0 = {\mathcal C}^{\widetilde L}$,
i.e., $\widetilde L$ is a pencil representation of the set
$\widetilde{\mathcal C} = {\mathcal C} - x^0$,
with $0$ in the interior of $\widetilde{\mathcal C}$.
By Lemma \ref{thm:monic} we may assume that $\widetilde L$ is a monic pencil.
The determinant
$\check{\tilde p}$ of $\widetilde L$
satisfies the $RZ_0$ property.  Thus the minimal defining
polynomial $\tilde p$ of $\widetilde{\mathcal C}$, which is the
shift of the minimal defining polynomial $p$ of ${\mathcal C}$
and which is contained in the $RZ_0$ polynomial
$\check{\tilde p}$ as a factor, satisfies $RZ_0$.  Shift $0$ back
to
$x^0$ to get that $p$ satisfies $RZ_{x^0}$.

The proof of the second assertion, for $m=2$, proceeds analogously. \hfill \qed

Notice a little more general than saying that the shifted pencil
$\widetilde L$ is monic, is saying that $\widetilde
L_0$ is PD.  This is equivalent to saying that
$L_0 + x^0_1 L_1 + \cdots + x^0_m L_m$ is PD,
and then, by the remark at the end of Section \ref{sec:necessity},
the interior of ${\mathcal C} = {\mathcal C}^L$ is exactly
$\{x : L_0 + x_1 L_1 + \cdots x_m L_m \mbox{ is PD}\}$.

\section{Geometrical  Viewpoint and
Examples of Sets  with no LMI Representation}
\label{sec:geom}

In this section we give a geometric characterization
of sets ${\mathcal C}$ with an LMI representation,
which is equivalent to Theorem \ref{thm:mainRZ}.

\subsection{Rigid Convexity}
\label{sec:rigidcon}


An algebraic interior  ${\mathcal C}$ of degree $d$ in
${\mathbb R}^m$ with minimal defining polynomial $p$
will be called {\bf rigidly convex} provided
for every point $x^0$ in the interior of ${\mathcal C}$
and a generic line
$\ell$ through $x^0$,
$\ell$ intersects the (affine) real algebraic hypersurface $p(x)=0$
in exactly $d$ points. { \footnote{
``Generic'' means that the exceptional lines are contained in a proper algebraic subvariety
(in ${\mathbb P}^{m-1}({\mathbb R})$); when $m=2$ this means simply all but a finite number
of lines.
One can replace here ``generic line'' by ``every line'' if one takes multiplicities into
account when counting the number of intersections, and also
counts the intersections at infinity, i.e., replaces the affine real
algebraic hypersurface $p(x_1,\dots,x_m)=0$ by the projective
real algebraic hypersurface $X_0^d
p(X_1/X_0,\dots,X_m/X_0)=0$.} }
We soon  show  that rigid convexity of $\cC_p$ is the same as
$p$ having the (shifted) \dRZ \ Property and as a consequence we obtain the following
{\bf geometric version of our main theorem.}

\begin{thm}
\label{thm:Geometric}
If ${\mathcal C}$  is a closed convex set in ${\mathbb R}^m$ with an
LMI representation, then ${\mathcal C}$ is rigidly
convex.  When $m=2$, the converse is true, namely, a rigidly
convex degree $d$ set has an LMI representation with
symmetric matrices $L_j \in {\mathbb R}^{d\times d}$.

\end{thm}

\noindent{\bf Proof}.

The theorem follows from Theorem \ref{thm:mainRZ}
(more precisely, its shifted version Corollary \ref{thm:shifted}) provided we can show
that {\it an algebraic interior $\cC$ with minimal defining polynomial $p$
is rigidly convex if and only if $p$ satisfies the shifted \dRZ \
condition.} We now set about to prove this equivalence.

Let $\ell: = \{x\in {\mathbb R}^m : x = \mu v + x^0,\ \mu \in {\mathbb R}\}$, $v\in {\mathbb R}^m$,
be a parameterization of a line through
$x^0 $ in  the interior of $ {\cC}$.
The points of intersection of the real algebraic hypersurface $p(x)=0$ with the line $\ell$
are parameterized
by exactly those $\mu \in {\mathbb R}$ at which
$$f{(\mu)}:=p(\mu v + x^0) =0,$$
and rigid convexity of $\cC_p$ says that (for a generic $\ell$)
there are exactly $d$ such distinct $\mu$.
However, the degree of $p$ is
$d$, so for a generic direction $v \in {\mathbb R}^m$, the degree of
the polynomial $f$ equals $d$
(the exceptional directions are given by $p_d(v)=0$,
where $p_d$ is the sum of degree $d$ terms in $p$).
Thus these $f$,
by the Fundamental Theorem of Algebra, have
exactly $d$ zeroes  $\mu \in \mathbb{C}$,
counting multiplicities.
Furthermore, for a generic direction $v$, all the zeroes are simple:
multiple zeroes correspond to common zeroes of $f$ and
$f' = \sum_{i=1}^m v_i \frac{\partial p}{\partial x_i}$,
that is to lines $\ell$ which pass through a point of intersection of
algebraic hypersurfaces $p(x)=0$ and
$\sum_{i=1}^m v_i \frac{\partial p}{\partial x_i} (x) = 0$;
but since $p$ is irreducible and
$\deg \sum_{i=1}^m v_i \frac{\partial p}{\partial x_i} < \deg p$,
this intersection is a proper algebraic subvariety of $p(x)=0$
hence of dimension at most $m-2$.
(This argument applies to both real and complex zeroes, by considering
real algebraic hypersurfaces and their complexifications, respectively.)
Thus rigid
convexity says precisely that all zeroes of $f$ are at real $\mu$,
for $f$ arising from a generic $v$.

To obtain that this is  equivalent to the $RZ_{x^0}$ Property for $p$
we give an elementary continuity argument which implies that all
$f$ have only real zeroes.
Pick
an exceptional direction $v$ (which might make $f$ have degree less than $d$):
if $f_n$ is a sequence of polynomials of degree $d$ in one variable with real coefficients
such that every $f_n$ has only real zeroes, and if the polynomials $f_n$ converge
coefficientwise to a polynomial $f$ (of degree less than or equal to $d$),
then $f$ has only real zeroes.


\hfill \qed

Notice that it follows from Theorem \ref{thm:Geometric} that for $m=2$ a rigidly convex set
is in fact convex, and (as in the Remark following Corollary \ref{thm:shifted}),
it is enough to require the defining property of rigid convexity to hold for lines
through a single point of ${\mathcal C}$.
Again, this holds for any dimension $m$ and can be proved directly using
the alternative topological characterization
of $RZ$ polynomials --- see Section \ref{sec:Homprops}.
We also indicate there a strong stability property of rigidly convex algebraic interiors,
and speculate about this stability being characteristic of rigid convexity.

\subsection{Examples}
\label{sec:examples}

\subsubsection{
Example 1.} \ \ The polynomial
$$p(x_1,x_2) =    x_1^3 - 3 x_2^2 x_1 - (x_1^2 + x_2^2 ) ^2
$$ has zero set shown in Figure \ref{fig:1L}.

The complement to $p=0$ in ${\mathbb R}^2$ consists of 4 components, 3 bounded
convex components where $p>0$ and an unbounded component where $p<0$.
Let us analyze one of the bounded components,
say the one in the right half plane,
${\mathcal C}$ is the closure of
$$\{ x: p(x) > 0, x_1 >0 \}.$$
Does ${\mathcal C}$ have an LMI representation?
To check this: fix a point $O$ inside ${\mathcal C}$, e.g. $O= (.7, 0)$.

\bigskip

\unitlength=0.3mm
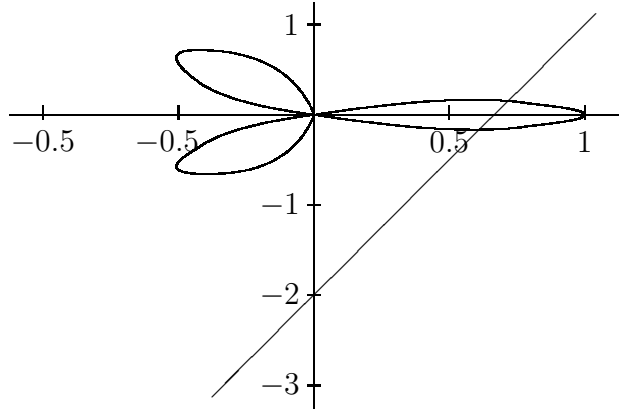
\begin{figure}[h]
\begin{picture}(400,200)(10,20)
\put(225,20){\line(0,1){180}}
\put(90,150){\line(1,0){270}}
\put(222,190){\line(1,0){06}}
\put(215,190){\makebox(0,0){{$1$}}}
\put(222,110){\line(1,0){06}}
\put(210,110){\makebox(0,0){{$-1$}}}
\put(222,70){\line(1,0){06}}
\put(210,70){\makebox(0,0){{$-2$}}}
\put(222,30){\line(1,0){06}}
\put(210,30){\makebox(0,0){{$-3$}}}
\put(105,148){\line(0,1){06}}
\put(105,138){\makebox(0,0){{$-0.5$}}}
\put(165,148){\line(0,1){06}}
\put(160,138){\makebox(0,0){{$-0.5$}}}
\put(285,148){\line(0,1){06}}
\put(285,138){\makebox(0,0){{$0.5$}}}
\put(345,148){\line(0,1){06}}
\put(345,138){\makebox(0,0){{$1$}}}
\qbezier(185,178)(215,174)(225,150)
\qbezier(173,165)(186,155)(225,150)
\qbezier (173,165)(150,182)(185,178)
\qbezier(185,124)(215,126)(225,150)
\qbezier(173,135)(186,145)(225,150)
\qbezier (173,135)(150,122)(185,124)
\qbezier(225,150)(290,160)(320,155)
\qbezier(225,150)(290,140)(320,145)
\qbezier(320,155)(370,150)(320,145)
\put(350,195){\line(-1,-1){170}}
\end{picture}
\caption{
$p =    x_1^3 - 3 x_2^2 x_1 - (x_1^2 + x_2^2 ) ^2$ with a line through $O=(.7,0)$
hitting $Z_p$ in only 2 points.
}
\label{fig:1L}
\end{figure}

By Theorem \ref{thm:Geometric}
almost every line $l$ thru $(0, 0.7)$ ( as in Figure \ref{fig:1L})
must intersect $p=0$ in 4 real points or the $RZ_O$
condition is violated.
We can see from the picture in $\bR2$ that
there is a continuum of real lines $\ell$ through $(0, 0.7)$
intersecting $p=0$
in exactly two real points.
Thus by Theorem \ref{thm:Geometric}  the set ${\mathcal C}$ does not
have an LMI
representation.
(Since $p$ is irreducible it is the minimum defining polynomial for
${\mathcal C}$.)

\subsubsection{
Example 2.}
$$p(x_1,x_2) = 1 -  x_1^4  - x_2^4  $$
Clearly, $\cC_p:= \{ x: p(x)\geq 0 \}$ has degree 4 but all
lines in $\bR2$ through it intersect
the set $p=0$
in exactly two places. Thus $\cC_p$ is not
rigidly convex.

\unitlength=0.20mm
\begin{figure}[h]
\begin{picture}(450,200)(10,20)
\put(225,10){\line(0,1){210}}
\put(75,115){\line(1,0){300}}
\put(222,215){\line(1,0){06}}
\put(210,215){\makebox(0,0){{$1$}}}
\put(222,165){\line(1,0){06}}
\put(210,165){\makebox(0,0){{$0.5$}}}
   \put(222,65){\line(1,0){06}}
\put(205,65){\makebox(0,0){{$-0.5$}}}
\put(222,15){\line(1,0){06}}
\put(205,15){\makebox(0,0){{$-1$}}}
\put(85,112){\line(0,1){06}}
\put(85,105){\makebox(0,0){{$-1$}}}
\put(155,112){\line(0,1){06}}
\put(155,105){\makebox(0,0){{$-0.5$}}}
\put(295,112){\line(0,1){06}}
\put(295,105){\makebox(0,0){{$0.5$}}}
\put(365,112){\line(0,1){06}}
\put(365,105){\makebox(0,0){{$1$}}}
\put(85,150){\line(0,1){10}}
\qbezier(85,160)(85,215)(165,215)
\put(165,215){\line(1,0){130}}
\put(365,150){\line(0,1){10}}
\qbezier(365,160)(365,215)(295,215)
\put(85,70){\line(0,1){10}}
\qbezier(85,70)(85,15)(165,15)
\put(165,15){\line(1,0){130}}
\put(365,70){\line(0,1){10}}
\qbezier(365,70)(365,15)(295,15)
\end{picture}
\caption{ $p(x_1,x_2) = 1 -  x_1^4  - x_2^4  $}
\label{fig:3}
\end{figure}
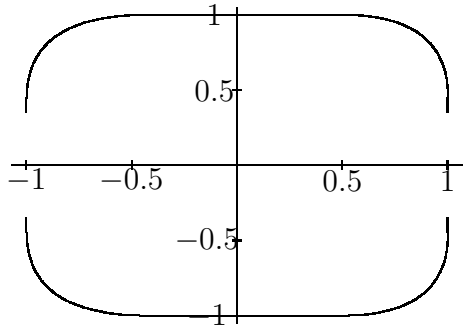


\section{Proof of the converse side of Theorem \ref{thm:mainRZ}}

A proof is an adaptation of the results of \cite{V93}.
Theorem 6.1 of \cite{V93} (see also \cite{Dub83})
says that  any polynomial $p$ on $\bR2$ satisfying
(\dRZ) (with $p(0)>0$) admits a representation (\ref{eq:detL}) with $L_0$, $L_1$ and $L_2$
complex hermitian matrices on ${\mathbb C}^d$, $d = \deg p$,
where $L_0$ is PD.
The possibility to take $L_0$, $L_1$ and $L_2$ real
symmetric follows easily from the arguments given in the last 4 paragraphs of
the Introduction (p.~456) of \cite{V93}.
It is assumed in \cite{V93} that $p$ is an irreducible polynomial, and the
projective closure of $p(x_1,x_2)=0$ is a smooth projective plane curve;
however the irreducibility assumption is irrelevant here since we may consider
one irreducible factor of $p$ at a time (and then form a direct sum of the
corresponding determinantal representations (\ref{eq:detL})),
and the smoothness assumption can be easily removed by considering {\em
maximal} determinantal representations as in \cite{BV96,BV99}.

This proof might lend itself to construction of the pencil $L$
via a construction of a basis for a linear system (or a Riemann Roch space)
on a plane algebraic curve.
\hfill\qed \vspace{3mm}

\noindent
Since the proof above is actually an outline which leaves
much for the reader to fill in,
we present an {\bf alternative direct proof}
based on the results of \cite[Section 5]{BV99}.
As before, we may assume without loss of generality
that $p$ is an irreducible polynomial;
we consider the complex affine algebraic curve
$\{(x_1,x_2)\in{\mathbb C}^2: p(x_1,x_2)=0\}$,
take its projective closure in the complex projective plane ${\mathbb P}^2({\mathbb C})$,
and let $M$ be the desingularizing Riemann surface.
$M$ is a compact real Riemann surface, i.e., a compact Riemann surface equipped with
an anti-holomorphic involution induced by
$(x_1,x_2) \longrightarrow (\bar x_1, \bar x_2)$. The set $M_{\mathbb R}$ of real points of $M$
is by definition the set of fixed points of the anti-holomorphic involution, and $x_1$
and $x_2$ are real meromorphic functions on $M$, i.e., meromorphic functions taking real values
at real points. Furthermore, it follows from the \dRZ \ condition that the real meromorphic
function $\frac{x_2}{x_1}$ on $M$ takes real values at real points only.
It follows that $M$ is a compact real Riemann surface of dividing type, i.e.,
$M \backslash M_{\mathbb R}$ consists of two connected components. We denote by $M_+$ and $M_-$
the components where $\Im\frac{x_2}{x_1} > 0$ and $\Im\frac{x_2}{x_1} < 0$, respectively.
We orient $M_{\mathbb R}$ so that $M_{\mathbb R} = \partial M_+$; then the real meromorphic
differential $-d\dfrac{x_1}{x_2} = \dfrac{x_1\,dx_2 - x_2\,dx_1}{x_2^2}$ is everywhere nonnegative
on $M_{\mathbb R}$.

Since both the assumptions and the conclusions of the theorem are invariant under a linear
change of coordinates $x_1$, $x_2$, we may assume without loss of generality that
$p(0,\mu)=0$ has $d$ distinct real roots $c_1, \ldots, c_d$; they correspond to $d$ non-singular
points $(0,c_1)$, \dots, $(0,c_d)$ on the affine algebraic curve $p(x_1,x_2)=0$, and we denote
by $q_1$, \dots, $q_d$ the corresponding points on $M$.
We may now write down explicit formulae for the $d \times d$ matrices $L_0$, $L_1$, $L_2$:
\begin{align}
& L_0 = I_d, \label{explicit0}\\
& [L_1]_{ij} = \begin{cases}
i=j: & \dfrac{\dfrac{dx_2}{dx_1}(q_i)}{c_i} \\
i \neq j: &
\left(\dfrac{1}{c_i}-\dfrac{1}{c_j}\right)
\dfrac{\theta\begin{bmatrix}a\\b\end{bmatrix}(\phi(q_j)-\phi(q_i))}
{\theta\begin{bmatrix}a\\b\end{bmatrix}(0) E(q_j,q_i)}
\dfrac{1}{\sqrt{d\dfrac{x_1}{x_2}}(q_i)\sqrt{d\dfrac{x_1}{x_2}}(q_j)}
\end{cases}, \label{explicit1}\\
& L_2 = \text{diag} \left(-\frac{1}{c_1}, \ldots, -\frac{1}{c_d}\right).\label{explicit2}
\end{align}
Here $\theta\begin{bmatrix}a\\b\end{bmatrix}(\cdot)$ is the theta function with characteristics
associated to the Jacobian variety of $M$, $\phi$ is the Abel-Jacobi map from $M$ into
the Jacobian variety, and $E(\cdot,\cdot)$ is the prime form on $M$; for details see
\cite[Section 2]{V93} or \cite[Section 4]{BV99}. Furthermore, we choose a symmetric homology
basis on $M$ as in \cite[Section 2]{V93} with a resulting period matrix $\Omega$, and
$b + \Omega a \in {\mathbb C}^g$ is a half period in the Jacobian with
$a_1, \ldots, a_{g+1-k} = 1$ and $b_{g+2-k}, \ldots, b_g = 0$ ($g$ is the genus of $M$ and
$k$ is the number of connected components of $M_{\mathbb R}$).

The fact that $\det(L_0+x_1L_1+x_2L_2)=p(x_1,x_2)$
(up to a constant positive factor) is exactly
\cite[Theorem 5.1]{BV99} (with a pair of meromorphic functions
$\lambda_1=1/x_1$, $\lambda_2=x_2/x_1$ on $M$ and $-x_1/x_2$ as a
local parameter at the poles of $\lambda_1$ and $\lambda_2$
(zeroes of $x_1$)). The fact that $L_1$ is a real symmetric matrix
follows immediately from the standard properties of theta
functions and \cite[Proposition 2.3]{V93}. Notice that $b + \Omega
a \in T_0$ hence $\theta\begin{bmatrix}a\\b\end{bmatrix}(0) \neq
0$, see \cite[Theorem 3.1 and Corollary 4.3]{V93}. \hfill\qed

\section{Topological Characterization of $RZ$ Polynomials}
\label{sec:oval}

This section gives an equivalent topological characterization
of zero sets of $RZ$ polynomials.
First
we show that $RZ$ polynomials of degree $d$
for which the sets of real zeros are compact
are exactly those whose zero sets
consist of a nest of $d/2$ ``ovaloids''
(having compact zero set implies $d$ even).
A result of similar definitiveness
holds when real zero sets are not compact,
but it requires projective space to neatly account for the effects of
the zero set passing thru infinity.
Thus we describe the general case only after describing
$p$ with compact real zero set, despite generating a little redundancy

Then we give some related properties of $RZ$ polynomials.

\subsection{$RZ$ Polynomials whose Sets of Real Zeroes are Compact}

%
%

We first define carefully what is meant by an ``ovaloid'' in  ${\mathbb R}^m$,
and later give an analogous definition in projective space, see for example
\cite{Viro86}.
We call $W \subset {\mathbb R}^m$ an {\em ovaloid in}  ${\mathbb R}^m$
if
$W$ is isotopic in  ${\mathbb R}^m$ to a sphere $S \subset {\mathbb R}^m$,
i.e., there is a homeomorphism $F$ of ${\mathbb R}^m$
with $F(S) = W$,
and furthermore $F$ is homotopic to
the identity,
i.e., there is a homeomorphism $H$ of
$[0,1] \times {\mathbb R}^m$ such that
$H_t=H|_{\{t\} \times {\mathbb R}^m}$ is a homeomorphism
of ${\mathbb R}^m $ for every $t$,
$H_0 = \operatorname{Id}_{{\mathbb R}^m}$,
and $H_1 = F$.
Notice that ${\mathbb R}^m \setminus S$ consists of two connected components
only one of which is contractible, hence the same is true
of ${\mathbb R}^m \setminus W$; we call the contractible component
the {\em interior} of the ovaloid $W$, and the
non-contractible component the {\em exterior}.
%

The following is a consequence of Theorem \ref{thm:oval} below.

\begin{cor}
\label{cor:oval}
Let $p$ be a polynomial
whose set $V$ of real zeros is compact
(then the degree $d$ of $p$ must be even, $d=2k$)
and which defines a smooth (affine) real algebraic hypersurface,
i.e., $p$ has no critical points on $V$.
Then $p$ satisfies $RZ_{x^0}$ with $p(x^0) \neq 0$ if and only if
$V$ is a disjoint union of $k$ ovaloids $W_1,\ldots,W_k$ in ${\mathbb R}^m$,
with $W_i$ contained in the interior of $W_{i+1}$, $i=1,\ldots,k-1$,
and $x^0$ lying in the interior of $W_1$.
Furthermore, in this case every  line $L$ through $x^0$ intersects
$V$ in $d$ distinct points.
\end{cor}

The corollary misses being as strong as possible
and in fact being sharp because of the requirement
that ovaloids be disjoint.
Intuitively, if
the real zero set of a
RZ polynomial is a  compact
singular (nonsmooth) hypersurface, then it
consists of nested ovaloids that touch at the
singular points, as in the figure below.
Notice that it follows from Property \ref{top-prop1} in Section \ref{sec:Homprops}
below that any $RZ$ polynomial of degree $d=2k$
with a compact real zero set is the limit
of $RZ$ polynomials of degree $d$ with compact real zero sets and which define
smooth real algebraic hypersurfaces (that is, the zero sets are 
disjoint unions of $k$ nested ovaloids).

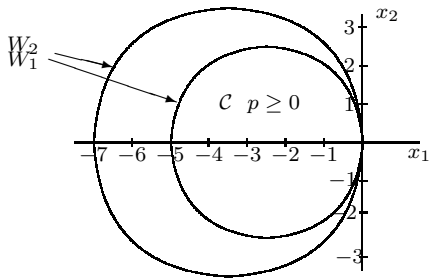
\unitlength=0.17mm
\begin{picture}(100,200)(10,20)
\put(355,0){\line(0,1){200}}
\put(130,100){\line(1,0){270}}

\qbezier(250,205)(350,200)(355,100)
\qbezier(250,-05)(350,0)(355,100)
\qbezier(250,-05)(150,0)(145,100)
\qbezier(250,205)(150,200)(145,100)

\qbezier(280,175)(350,170)(355,100)
\qbezier(280,25)(350,30)(355,100)
\qbezier(280,25)(210,30)(205,100)
\qbezier(280,175)(210,170)(205,100)

\put(352,190){\line(1,0){06}}
\put(345,190){\makebox(0,0){{\tiny{$3$}}}}
\put(353,160){\line(1,0){06}}
\put(345,160){\makebox(0,0){{\tiny{$2$}}}}
\put(353,130){\line(1,0){06}}
\put(345,130){\makebox(0,0){{\tiny{$1$}}} }
\put(353,70){\line(1,0){06}}
\put(340,70){\makebox(0,0){{\tiny{$-1$}}}}
\put(353,45){\line(1,0){06}}
\put(340,45){\makebox(0,0){{\tiny{$-2$}}}}
\put(353,10){\line(1,0){06}}
\put(345,10){\makebox(0,0){{\tiny{$-3$}}}}

\put(325,97){\line(0,1){06}}
\put(325,90){\makebox(0,0){{\tiny{$-1$}}}}
\put(295,97){\line(0,1){06}}
\put(295,90){\makebox(0,0){{\tiny{$-2$}}}}
\put(265,97){\line(0,1){06}}
\put(265,90){\makebox(0,0){{\tiny{$-3$}}}}
\put(235,97){\line(0,1){06}}
\put(235,90){\makebox(0,0){{\tiny{$-4$}}}}
\put(205,90){\makebox(0,0){{\tiny{$-5$}}}}
\put(175,97){\line(0,1){06}}
\put(175,90){\makebox(0,0){{\tiny{$-6$}}}}
\put(145,90){\makebox(0,0){{\tiny{$-7$}}}}

\put(90,176){\makebox(0,0){{\tiny{$W_2$}}}}
\put(90,164){\makebox(0,0){{\tiny{$W_1$}}}}

\put(275,130){\makebox(0,0){{\tiny{${\mathcal C}\ \ p \ge 0$}}}}
\put(115,170){\vector(4,-1){45}}
\put(110,165){\vector(3,-1){100}}

\put(375,200){\makebox(0,0){{\tiny{$x_2$}}}}
\put(400,90){\makebox(0,0){{\tiny{$x_1$}}}}
\end{picture}


$$p(x_1, x_2) =   ( x_1^2+ x_2^2 )(x_1^2+   x_2^2 + 12 x_1 - 1)+36 x_1^2 \geq 0
$$


\subsection{$RZ$ Polynomials with a Non-Compact Set of Real Zeroes}

We shall
need to consider the behavior of the polynomial at infinity,
so we introduce the homogenization $P$ of $p$,
$$P(X_0,X_1,\ldots,X_m)=X_0^d \; p(X_1/X_0,\ldots,X_m/X_0),$$
where $d = \deg p$.
Notice
that $P$ is not divisible by $X_0$, and that conversely given
any homogeneous polynomial $P$ of degree $d$ in $m+1$ variables
which is not divisible by $X_0$ we may recover a polynomial $p$ of degree $d$
in $m$ variables so that $P$ is the homogenization of $p$ by
$$p(x_1,\ldots,x_m)=P(1,x_1,\ldots,x_m).$$
We identify as usual the $m$ dimensional real projective
space ${\mathbb P}^m({\mathbb R})$ with
the union of ${\mathbb R}^m$ and the hyperplane
at infinity $X_0=0$,
so that the affine coordinates $x$ and the
projective coordinates $X$ are related by
$x_1=X_1/X_0$, \dots, $x_m=X_m/X_0$;
the {\it projective real algebraic hypersurface} $P(X)=0$ in
${\mathbb P}^m({\mathbb R})$ is called the {\it projective closure}
of
the affine real algebraic hypersurface $p(x)=0$ in ${\mathbb R}^m$.
We shall say that $p$ defines a {\it smooth projective real algebraic hypersurface}
if $P$ has no critical points on $P(X)=0$.
We shall often abuse the notation and denote a point in ${\mathbb P}^m({\mathbb R})$
and a $(m+1)$-tuple of its projective coordinates by the same letter.

We call $W \subset {\mathbb P}^m({\mathbb R})$ an {\em ovaloid} if
$W$ is isotopic in ${\mathbb P}^m({\mathbb R})$
to a sphere $S \subset {\mathbb R}^m \subset {\mathbb P}^m({\mathbb R})$,
i.e., there is a homeomorphism $F$ of ${\mathbb P}^m({\mathbb R})$
with $F(S) = W$, and furthermore $F$ is homotopic to
the identity, i.e., there is a homeomorphism $H$ of
$[0,1] \times {\mathbb P}^m({\mathbb R})$ such that
$H_t=H|_{\{t\} \times {\mathbb P}^m({\mathbb R})}$ is a homeomorphism
of ${\mathbb P}^m({\mathbb R})$ for every $t$,
$H_0 = \operatorname{Id}_{{\mathbb P}^m({\mathbb R})}$,
and $H_1 = F$.
Notice that ${\mathbb P}^m({\mathbb R}) \setminus S$ consists of two connected components
only one of which is contractible, hence the same is true
of ${\mathbb P}^m({\mathbb R}) \setminus W$; we call the contractible component
the {\em interior} of the ovaloid $W$, and the non-contractible component the {\em exterior}.

We call $W \subset {\mathbb P}^m({\mathbb R})$ a {\em pseudo-hyperplane}
if $W$ is isotopic in ${\mathbb P}^m({\mathbb R})$ to a (projective) hyperplane
$H \subset {\mathbb P}^m({\mathbb R})$.

\begin{theorem}
\label{thm:oval}
Let $p$ be a polynomial of degree $d$ in $m$ variables
which defines a smooth projective real algebraic hypersurface
$V$ in ${\mathbb P}^m({\mathbb R})$.
Then $p$ satisfies $RZ_{x^0}$ with $p(x^0) \neq 0$ if and only if
\begin{itemize}
\item[a.]
if $d=2k$ is even, $V$ is a disjoint union of $k$ ovaloids $W_1,\ldots,W_k$,
with $W_i$ contained in the interior of $W_{i+1}$, $i=1,\ldots,k-1$,
and $x^0$ lying in the interior of $W_1$;
\item[b.]
if $d=2k+1$ is odd, $V$ is a disjoint union of $k$ ovaloids $W_1,\ldots,W_k$,
with $W_i$ contained in the interior of $W_{i+1}$, $i=1,\ldots,k-1$,
and $x^0$ lying in the interior of $W_1$, and a pseudo-hyperplane $W_{k+1}$
contained in the exterior of $W_k$.
\end{itemize}
Furthermore, in this case every (projective) line $L$ through $x^0$ intersects
$V$ in $d$ distinct points.
\end{theorem}

The proof  appears in  Section \ref{sec:oval proofs}.

\subsection{Topological Properties of $RZ$ Polynomials}
\label{sec:Homprops}

We list now some properties of $RZ$ polynomials, mostly of a topological
nature, which are closely related to Theorem \ref{thm:oval}.

\noindent
{\bf RZ Topological Properties}

\begin{enumerate}

\item
\label{top-prop1}
If two polynomials $p_0, p_1$ of degree $d$ in $m$ variables satisfy \dRZ \
and $p_0(0), p_1(0) > 0$, then
there is a continuous path through \dRZ \ polynomials of degree $d$ in $m$ variables
from $p_0$ to $p_1$
with each polynomial $p$ in the path satisfying $p(0)>0$.
Furthermore, we may assume that except for possibly
the endpoints $p_0, p_1$, the polynomial $p$
defines a smooth projective real algebraic hypersurface.


\item
\label{top-prop3}
If the polynomial $p$ satisfies $RZ_{\tilde x^0}$ for some point $\tilde
x^0$ and
$p(\tilde x^0) > 0$, then $p$ satisfies $RZ_{x^0}$ for every interior point $x^0$
of ${\mathcal C}_p(\tilde x^0)$.

\item
\label{top-prop4}
If $p$ is a \dRZ \ polynomial with $p(0)>0$ then the algebraic interior
${\mathcal C}_p$ is convex
(so a rigidly convex algebraic interior is in fact convex).

\item
\label{top-prop5}
Suppose that
there is a continuous path through polynomials of degree $d$
in $m$ variables from a polynomial
$p_0$ to a polynomial $p_1$, such that
each polynomial $p$ in the path
defines a smooth projective real algebraic hypersurface and $p(0) > 0$.
Then if $p_0$ is a \dRZ \ polynomial, so is $p_1$.

\item
\label{top-prop6}
Suppose that
there is a continuous path through polynomials of degree $d$
in $m$ variables from a polynomial
$p_0$ to a polynomial $p_1$, such that
each polynomial $p$ in the path including $p_0$ and $p_1$
defines a smooth projective real algebraic hypersurface.
Assume that $p_0$ satisfies $RZ_{x^0}$ for some point $x^0$ and $p(x^0) > 0$;
then $p_1$ satisfies $RZ_{x^1}$ for some point $x^1$ and $p(x^1) > 0$.
\end{enumerate}

The proof  appears in  Section \ref{sec:rzproofs}.

\section{Generalizations}
\label{sec:fragen}

We now conjecture structure which extends that actually proved in this paper.
\subsection{{\bf $m>2$}}
It is a natural question to which extent the converse side of
our main results Theorem \ref{thm:mainRZ} and Theorem \ref{thm:Geometric}
extends to dimension $m>2$. A simple count of parameters shows that given
a polynomial $p(x_1,\dots,x_m)$ in $m > 2$ variables, in general there do not
exist $d \times d$ matrices $L_0,\dots,L_m$ such that
$p(x) = \det (L_0 + x_1 L_1 + \cdots + x_m L_m)$
(here the polynomial has complex coefficients and the matrices are likewise complex).
Thus we cannot expect a rigidly convex degree $d$ set in ${\mathbb R}^m$, $m>2$,
to admit an LMI representation with $d \times d$ symmetric matrices.
However the count of parameters does not preclude the existence of a pencil representation
with matrices of a size larger than the degree.
In fact we conjecture that
\begin{quote}
{\it
A convex set $\mathcal C $ in ${\mathbb R}^m$ has an LMI representation if and
only if it is rigidly convex.
}
\end{quote}

\subsection{LMI Lifts}
Another practically important question (in the spirit of
Nestrov and Nimmerovski \cite{NN94} )
is

\noindent
{\bf
Question}
{\it Find a test which insures that a
convex set $\mathcal C$ in ${\mathbb R}^m$ lifts to
some LMI representable set ${\ctC}$ in a possibly
bigger space ${\mathbb R}^{m+k}$.
That is,
find neccesary and sufficient properties
on a given set
 $\mathcal C$ in ${\mathbb R}^m$
which insure that there exists
an LMI representable set ${\ctC}$ in ${\mathbb R}^{m+k}$
whose projection onto ${\mathbb R}^m$ equals the set $\mathcal C$?
}

No apriori restriction on such a $\mathcal C$ is appearent,
either from theoretical
considerations or from numerical experiments
(done informally by Pablo Parrilo).
Of course, $\mathcal C$ has to be a semi-algebraic set with a connected interior,
which is equal to the closure of its interior. It is not clear whether $\mathcal C$
should be an algebraic interior as we defined it.

\subsection{Convexity under Deformation is Ephemeral}

Properties \ref{top-prop5} and \ref{top-prop6} in Section  \ref{sec:Homprops}
imply that the convexity of a rigidly convex algebraic interior
${\mathcal C}$ is stable under deformation of the minimal defining polynomial $p$
as long as the projective real algebraic hypersurface defined by $p$
(the projective closure in ${\mathbb P}^m({\mathbb R})$ of the Zariski
closure of $\partial{\mathcal C}$ in ${\mathbb R}^m$) remains smooth.
We conjecture that this property characterizes rigidly convex algebraic interiors;
more precisely, we make the following conjecture:
\begin{quote}
{\it
Let ${\mathcal C}_0 = {\mathcal C}_{p_0}$ be a convex algebraic interior of
degree $d$ in ${\mathbb R}^m$ which is not rigidly convex and such
that the polynomial $p_0$ defines a smooth projective real algebraic hypersurface.
Then there exists a continuous path through polynomials of degree $d$
in $m$ variables from a polynomial
$p_0$ to a polynomial $p_1$, such that
each polynomial $p$ in the path
defines a smooth projective real algebraic hypersurface and $p(0) > 0$,
and the algebraic interior ${\mathcal C}_{p_1}$ is not convex.
}
\end{quote}

%

For more information about various connected components of
the set of polynomials of degree $d$
in $m$ variables that define a smooth projective real algebraic hypersurface,
see \cite{Viro86}.


\section{Appendix: Proofs of Topological Properties}

We first establish Property \ref{top-prop1} of Section \ref{sec:Homprops},
and use it to prove Theorem \ref{thm:oval}.
Then we establish Properties \ref{top-prop3}--\ref{top-prop6}.

\label{sec:rzproofs}

{\bf Proof of Property \ref{top-prop1}.}
This crucial arcwise connectivity property follows 
from the results of Nuij \cite{Nuij68} on hyperbolic polynomials.
Recall that a homogeneous polynomial $P$ in 
$m+1$ variables is said to be hyperbolic
with respect to $X^0=(X^0_0,X^0_1,\ldots,X^0_m)\in{\mathbb R}^{m+1}$ if
$P(X^0) \neq 0$ and for each $X \in {\mathbb R}^{m+1}$ all the zeroes $\lambda$
of $P(X + \lambda X^0)$ are real. It is called strictly hyperbolic if all the zeroes
are furthermore simple for each $X \in {\mathbb R}^{m+1}$ (not a multiple of
$X^0$), and it is called normalized if $P(X^0)=1$.

It is immediate that a polynomial $p$ in $m$ variables satisfies $RZ_{x^0}$
and $p(x^0) \neq 0$
for $x^0=(x^0_1,\ldots,x^0_m)$ if and only if its homogenization $P$
is hyperbolic with respect to $X^0=(1,x^0_1,\ldots,x^0_m)$.
Notice that if $P$ is strictly hyperbolic, then necessarily
$p$ defines a smooth projective real algebraic hypersurface $V$. \footnote{
It will follow from Theorem \ref{thm:oval} that the converse is also true ---
if $p$ defines a smooth projective real algebraic hypersurface $V$,
then $P$ is strictly hyperbolic.}
For (as in the proof of Theorem \ref{thm:Geometric}) the zeroes of $P(X + \lambda X^0)$
correspond to intersections of the line through $X$ and $X^0$ with $V$;
if $X$ were a singular point of $V$, the intersection multiplicity at $X$
of any line through $X$ with $V$ would be greater than 1 and the corresponding
zero of $P(X + \lambda X^0)$ could not be simple.

Property \ref{top-prop1} now follows from the following results of \cite{Nuij68}
(where the polynomials are of degree $d$ and are hyperbolic with respect
to a fixed $X^0$):
(a) the set of strictly hyperbolic homogeneous polynomials is open;
(b) every hyperbolic homogeneous polynomial is the end point of an open continuous
path of strictly hyperbolic homogeneous polynomials;
(c) the set of normalized strictly hyperbolic homogeneous polynomials is
arcwise connected.
Notice that since $X^0=(1,x^0_1,\ldots,x^0_m)$,
homogeneous polynomials which are divisible by $X_0$
form an affine subspace of codimension at least 2 in the affine space
of homogeneous polynomials $P$ with $P(X^0)=1$;
so the set of normalized strictly hyperbolic homogeneous polynomials
which are not divisible by $X_0$ is still arcwise connected.
\qed

{\bf Proof of Theorem \ref{thm:oval}, first direction.}
Assume that $p_0$ is a $RZ$ polynomial of degree $d$ in $m$ variables with $p_0(0)>0$
which defines a smooth projective real algebraic hypersurface $V_0$ in
${\mathbb P}^m({\mathbb R})$.
Let
$$p_1(x_1,\ldots,x_m) = (r_1^2 - x_1^2 - \cdots - x_m^2) \cdots
                        (r_k^2 - x_1^2 - \cdots - x_m^2)$$
if $d=2k$ is even, and
$$p_1(x_1,\ldots,x_m) = (r_1^2 - x_1^2 - \cdots - x_m^2) \cdots
                        (r_k^2 - x_1^2 - \cdots - x_m^2) (r_{k+1} - x_1)$$
if $d=2k+1$ is odd, for some $0 < r_1 < \cdots < r_k < r_{k+1}$.
By Property \ref{top-prop1}, there exists a continuous path $\{p_t\}_{0 \leq t \leq 1}$
from $p_0$ to $p_1$
through $RZ$ polynomials of degree $d$ in $m$ variables such that $p_t(0)>0$
and $p_t$ defines a smooth projective real algebraic hypersurface $V_t$ in
${\mathbb P}^m({\mathbb R})$ for each $t$.
It follows from the Thom Isotopy Theorem (see \cite[Theorem 14.1.1]{BCR98} or
\cite[Theorem 20.2]{AbRo67}) that there exists a homeomorphism
$H$ of $[0,1] \times {\mathbb P}^m({\mathbb R})$ such that
$H_t = H|_{\{t\} \times {\mathbb P}^m({\mathbb R})}$ is a homeomorphism
of ${\mathbb P}^m({\mathbb R})$ for every $t$,
$H_0 = \operatorname{Id}_{{\mathbb P}^m({\mathbb R})}$,
and $H_t(V_0) = V_t$.
Therefore $V_0$ is isotopic to $V_1$,
and $V_1$ by construction
is a disjoint union of $k$ nested ovaloids if $d=2k$
and a disjoint union of $k$ nested ovaloids and a pseudo-hyperplane if $d=2k+1$,
thus $V_0$ is also,
as is required.

It remains only to show that $0$ is contained in the interior of the innermost
ovaloid. Let $x^t = H_t^{-1}(0)$; since $p_t(0) > 0$ for all $t$,
$0 \not\in V_t$ and hence $x^t \not\in V_0$.
Thus $x^t$ is a continuous path in ${\mathbb P}^m({\mathbb R}) \setminus V_0$.
Since $0$ obviously belongs to the interior of the innermost sphere in $V_1$,
it follows that $x^1$ belongs to the interior of the innermost
ovaloid in $V_0$. Therefore so does $x^0=0$. \qed

{\bf Proof of Theorem \ref{thm:oval}, second direction.}
Assume now conversely that $p$ is a polynomial of degree $d$ in $m$ variables
that defines a smooth projective real algebraic hypersurface $V$
in ${\mathbb P}^m({\mathbb R})$ which is a disjoint union of $k$ nested ovaloids
$W_1,\ldots,W_k$ if $d=2k$
and a disjoint union of $k$ nested ovaloids $W_1,\ldots,W_k$
and a pseudo-hyperplane $W_{k+1}$ if $d=2k+1$. Assume that $X^0$ is in the interior of $W_1$,
and let $L$ be a (projective) line in ${\mathbb P}^m({\mathbb R})$ through $X^0$.

It is obvious that there exists a pseudo-hyperplane $W$ contained in the exterior
of the outermost ovaloid $W_k$ (just take the image of a hyperplane
contained in the exterior of a sphere $S$ in ${\mathbb P}^m({\mathbb R})$).
Since the ${\mathbb Z}/2$ fundamental class of $W$ coincides with that of
a hyperplane, it follows that the ${\mathbb Z}/2$ intersection number
of $L$ and $W$ equals 1, hence $L \cap W \neq \emptyset$. Let $X \in L \cap W$;
in case $d=2k+1$ is odd we take $W=W_{k+1}$ and denote $X^{k+1}=X$.

The points $X^0$ and $X$ decompose the projective line $L$ into two closed intervals
$L_+$ and $L_-$. Each of $L_\pm$ is a continuous path between $X^0$ which lies
in the interior of $W_1$ and $X$ which lies in the exterior of $W_k$; therefore
$L_\pm$ intersects each of $W_1,\ldots,W_k$ in at least one point $X^1_\pm,\ldots,X^k_\pm$.
Adding to these $X^{k+1}$ in case $d=2k+1$ is odd, we obtain $d$ distinct points
of intersection of $L$ with $V$; since the defining polynomial $p$ of $V$ has degree $d$,
$L \cap V$ can contain at most $d$ points, so we conclude that $L \cap V$ contains
exactly $d$ distinct points $X^1_\pm,\ldots,X^k_\pm$ and (in case $d=2k+1$ is odd) $X^{k+1}$.
This completes the proof.
Notice that it follows that on $L_\pm$ we have the following picture:
the interval $[X^0,X^1_\pm)$ lies in the interior of $W_1$,
the interval $(X^1_\pm,X^2_\pm)$ lies in the interior of $W_2$
intersection the exterior of $W_1$, \dots,
the interval $(X^{k-1}_\pm,X^k_\pm)$ lies in the interior of $W_k$
intersection the exterior of $W_{k-1}$, and
the interval $(X^k_\pm,X]$ lies in the exterior of $W_k$. \qed

{\bf Proof of Property \ref{top-prop3}.}
By Property \ref{top-prop1}, we may assume, using approximation,
that the polynomial $p$ defines a smooth projective real algebraic hypersurface.
The conclusion now follows immediately from Theorem \ref{thm:oval}. \qed

{\bf Proof of Property \ref{top-prop4}.}
We may assume, again,
that the polynomial $p$ defines a smooth projective real algebraic hypersurface $V$
in ${\mathbb P}^m({\mathbb R})$.
By Theorem \ref{thm:oval}, $V$ is a disjoint union of $k$ nested ovaloids
$W_1,\ldots,W_k$ if $d=2k$
and a disjoint union of $k$ nested ovaloids $W_1,\ldots,W_k$
and a pseudo-hyperplane $W_{k+1}$ if $d=2k+1$;
$X^0=(1,0,\ldots,0)$ (we are using here projective coordinates)
lies in the interior ${\mathcal I}$ of $W_1$.

It is obvious that ${\mathcal C}_p$ is the closure of
the connected component of $X^0$ in
${\mathcal I} \setminus {\mathcal I} \cap H_\infty$
where $H_\infty = \{ X_0 = 0 \}$ is the hyperplane at infinity.
We shall show that either ${\mathcal I} \cap H_\infty = \emptyset$ in which case
${\mathcal I} \subset {\mathbb R}^m$ is convex,
or ${\mathcal I} \setminus {\mathcal I} \cap H_\infty$ consists
of two connected convex components in ${\mathbb R}^m$.

Before proceeding, let us notice that there exists a hyperplane $H$
in ${\mathbb P}^m({\mathbb R})$ such that $H \cap {\mathcal I} = \emptyset$.
We may take $H$ to be the (projectivized) tangent hyperplane to $V$ at a point
$Y^0 \in W_1$. If $H$ contained a point $Y$ in ${\mathcal I}$, then
by Theorem \ref{thm:oval} the line through $Y$ and $Y^0$ would intersect $V$
in $d$ distinct points; but since this line is contained in $H$ it has
intersection multiplicity at least 2 with $V$ at $Y^0$, so it cannot
intersect $V$ in more than $d-1$ distinct points.

Let now $Y^1,Y^2 \in {\mathcal I} \setminus {\mathcal I} \cap H_\infty$.
Let $L$ be the projective line through $Y^1$ and $Y^2$ and let
$X^\infty = L \cap H_\infty$. The two points $Y^1$ and $Y^2$ decompose
the projective line $L$ into two closed intervals; the line segment
between $Y^1$ and $Y^2$ in ${\mathbb R}^m$ is the interval that does
not contain the point $X^\infty$. Notice that by the last sentence
in the proof of Theorem \ref{thm:oval}, exactly one of the two intervals is contained
in ${\mathcal I}$; let us denote it by $L_{\mathcal I}$.

Assume that ${\mathcal I} \cap H_\infty = \emptyset$.
Then $X^\infty \not\in {\mathcal I}$ and therefore
$X^\infty \not\in L_{\mathcal I}$.
It follows that $L_{\mathcal I}$ is the line segment
between $Y^1$ and $Y^2$ in ${\mathbb R}^m$,
and ${\mathcal I} \subset {\mathbb R}^m$ is convex.

If ${\mathcal I} \cap H_\infty \neq \emptyset$,
let $a_0 X_0 + a_1 X_1 + \cdots + a_m X_m = 0$ be a hyperplane in
${\mathbb P}^m({\mathbb R})$ that does not intersect ${\mathcal I}$.
Let ${\mathcal I}_+$ and ${\mathcal I}_-$ be the  nonempty open subsets of
${\mathcal I} \setminus {\mathcal I} \cap H_\infty$
where $X_0 / (a_0 X_0 + a_1 X_1 + \cdots + a_m X_m) > 0$
and $X_0 / (a_0 X_0 + a_1 X_1 + \cdots + a_m X_m) < 0$, respectively.
We shall show that ${\mathcal I}_+$ and ${\mathcal I}_-$ are convex subsets
of ${\mathbb R}^m$, thereby completing the proof.
It is enough to show that if the line segment
between $Y^1$ and $Y^2$ in ${\mathbb R}^m$ is not contained in ${\mathcal I}$
and say $Y^1 \in {\mathcal I}_+$, then necessarily $Y^2 \in {\mathcal I}_-$.

Assume then that the line segment
between $Y^1$ and $Y^2$ in ${\mathbb R}^m$ is not contained in ${\mathcal I}$.
This means that this line segment is not $L_{\mathcal I}$,
so $X^\infty \in L_{\mathcal I}$.
It is now immediate that the function
$X_0 / (a_0 X_0 + a_1 X_1 + \cdots + a_m X_m)$ on $L_{\mathcal I}$
changes its sign at $X^\infty$. \qed

{\bf Proof of Property \ref{top-prop5}.}
This follows immediately from the facts that a limit of $RZ$ polynomials is a $RZ$
polynomial (see the last sentence in the proof of Theorem \ref{thm:Geometric}),
that a $RZ$ polynomial $p$, $p(0)>0$, defines a smooth projective real algebraic hypersurface
if and only if its homogenization $P$ is strictly hyperbolic
(with respect to $X^0=(1,0,\ldots,0)$),
and that the set of strictly hyperbolic homogeneous polynomials is open \cite{Nuij68}. \qed

(We could also prove Property \ref{top-prop5}
by using Theorem \ref{thm:oval} and an isotopy argument, like the proof
of Property \ref{top-prop6} below.)

{\bf Proof of Property \ref{top-prop6}.}
By Thom Isotopy Theorem, smooth projective real algebraic hypersurfaces
$V_0$ and $V_1$ in ${\mathbb P}^m({\mathbb R})$ defined by the polynomials
$p_0$ and $p_1$ respectively are isotopic; the result now follows from Theorem \ref{thm:oval}. \qed

\label{sec:oval proofs}

\section{ Thanks }
We are grateful to Bernd Sturmfels
for suggesting  this problem
and to Pablo Parrilo for discussions bearing
on possible extensions of these results.

Thanks for Helton's support  are due to the NSF, DARPA, and the
Ford Motor company.

Thanks for Vinnikov's support are due to the Israel Science Foundation.


\providecommand{\bysame}{\leavevmode\hbox to3em{\hrulefill}\thinspace}

\newpage
\vspace{-.2in}
 \centerline{{\bf  NOT FOR PUBLICATION}}
 \vspace{-.2in}
\tableofcontents

\end{document}